\documentclass{amsart}
\makeatletter
\usepackage{geometry}
\newgeometry{margin=1.2in}
\usepackage{graphicx}
\usepackage{tikz-cd}
\usepackage{amsmath,amssymb}
\usepackage{mathtools}
\usepackage{hyperref}
\numberwithin{equation}{section}
\DeclareMathOperator{\Z}{\mathbb{Z}}

\DeclareMathOperator{\bl}{B\ell}
\DeclareMathOperator{\id}{id}
\DeclareMathOperator{\pp}{\mathbb{P}}
\DeclareMathOperator{\oo}{\mathcal{O}}
\DeclareMathOperator{\hilb}{Hilb}
\DeclareMathOperator{\sym}{Sym}
\DeclareMathOperator{\res}{res}

\newcommand{\wt}{\widetilde}
\newcommand{\ra}{\rightarrow}
\newcommand{\cal}{\mathcal}
\newtheorem{thm}{Theorem}[section]
\newtheorem{prop}[thm]{Proposition}
\newtheorem{cor}[thm]{Corollary}
\newtheorem{lm}[thm]{Lemma}

\theoremstyle{definition}
\newtheorem{defi}[thm]{Definition}
\theoremstyle{remark}
\newtheorem{rmk}[thm]{Remark}
\theoremstyle{theorem}
\newtheorem{introthm}{Theorem}

\title[The Chow ring of the Hilbert cube]{The Chow ring of the Hilbert cube}
\author[Ian Selvaggi]{Ian Selvaggi}
\address{Scuola Internazionale Superiore di Studi Avanzati (SISSA), Italy} %
\email{iselvagg@sissa.it}
\date{\today}
\begin{document}
\begin{abstract}
    Given a smooth projective variety $X$ over an algebraically closed field $k$, we compute the Chow ring of the Hilbert scheme of three points on $X$, $\hilb^3(X)$, as an algebra with generators and relations over the Chow ring of $X\times\sym^2(X)$. If in addition the characteristic of $k$ is zero, we extend the computation to the quasi-projective case.
\end{abstract}
\maketitle
\tableofcontents
\section{Introduction}
Hilbert and Quot schemes, originally introduced by Grothendieck in his foundational work \cite{grothendieck1957fondements}, are among the most studied moduli spaces in algebraic geometry. While tackling the study of their geometry in full generality is completely out of reach with the current techniques, in some specific contexts such a task becomes slightly more manageable. This is the case of Hilbert schemes of points, namely those moduli spaces parametrizing zero dimensional subschemes of a fixed length within a given variety $X$. 

Let us fix a ground field $k$ and a smooth quasi-projective $k$-variety $X$. The Hilbert scheme of $m$ points over $X$ is by definition the (quasi-projective) scheme representing the functor
\begin{align*}
    \hilb^m(X):(\operatorname{Sch}/k)^{op}&\ra(\operatorname{Set})\\
    T\mapsto\{\oo_{X_T}\twoheadrightarrow&\mathcal{Q}\,|\,\mathcal{Q}\hbox{ is }T\hbox{-flat},\,\chi(\mathcal{Q}|_{X_t})=m\hbox{ for all } t\in T\hbox{ geometric points}\},
\end{align*}
where $\chi(-)$ denotes the Euler characteristic. The behaviour of these spaces is increasingly complex as the dimension of $X$ and the number $m$ increase. For example:
\begin{itemize}
    \item when $X$ is a projective curve $\hilb^m(X)\cong\sym^m(X)$ for all $m$, and it is smooth if $X$ is,
    \item $\hilb^m(X)$ is smooth if and only if $\dim(X)\leq2$ by \cite{fogarty1968algebraic} or $m\leq 3$,
    \item $m=4$, $\dim (X)=3$ - $\hilb^4(\mathbb{A}^3)$ is already singular,
    \item $m$ and $\dim (X)$ arbitrary - \emph{hic sunt leones} (see \cite{jelisiejew2020pathologies}).
\end{itemize}
If understanding the geometry of Hilbert schemes is hard, computing their invariants, in general, is no easier job. For surfaces, (co)homological and motivic structures, and Chow rings are quite understood by \cite{ellingsrud1987homology}, \cite{gottsche1990betti}, \cite{nakajima1997heisenberg}, \cite{DECATALDO2002824}. In higher dimension instead, the overall situation is still very involved. On the motivic side, there has been some interest lately for $\hilb^m(\mathbb{A}^d)$ (see \cite{zhan2022punctual} and \cite{graffeo2024motive}), while some partial results are known for homology groups of desingularizations \cite{katz1994desingularization}.

If we fix $m=3$, as remarked above $\hilb^3(X)$ is known to be smooth and projective if $X$ is. With such hypotheses, Fantechi and Göttsche managed to compute the cohomology ring of $\hilb^3(X)$ as an algebra with generators and relations over the cohomology ring of $X$ \cite{Gottsche1993}. The purpose of this article is to extend their computation to the Chow ring $A^*(\hilb^3(X))$. If $\operatorname{char}(k)=0$, we are able to carry the result to the quasi-projective setting.

The strategy is very similar to the one in \cite{Gottsche1993}, and will be outlined in the next section. The only difficulty arises from a key property (or better, a lack of) of Chow groups: there is no Künneth formula for $A^*(X)$. Without further assumptions, the best we can do is to describe $A^*(\hilb^3(X))$ as an algebra over the ring $A^*(X\times\sym^2(X))$ as follows.
\begin{introthm}[Theorem \ref{thmchow3}]
    The Chow ring of $\hilb^3(X)$ with coefficients in $\Z\left[\frac{1}{6}\right]$ is a subring of an algebra over $A^*(X\times\sym^2(X))$ with two generators and relations \eqref{rel0}, \eqref{rel1}, \eqref{rel2}, \eqref{rel3}.
\end{introthm}
In the case where $X$ admits a Chow--Künneth decomposition, i.e. when $$A^*(X\times X)\cong A^*(X)\otimes A^*(X),$$
we recover the description of $A^*(\hilb^3(X))$ given in \cite{Gottsche1993}. However, it should be also remarked that with such hypotheses one can say much more on $A^*(\hilb^3(X))$ - see \cite{shen2016motive}.\\\\
\textbf{Notation and conventions.}
We follow the conventions in \cite{fulton2013intersection}. We work over an algebraically closed field $k$. For a $k$-variety $X$, we denote by $A_d(X)$ its $d$-th Chow group, that is the group of $k$-cycles modulo rational equivalence. If $X$ is smooth, set $A^d(X):=A_{\dim(X)-d}(X)$ and define the Chow ring as $A^*(X):=\bigoplus_{k=0}^{\dim(X)}A^k(X)$, with the ring structure given by intersection product. With some abuse of notation, for a map $f:X\ra Y$ of smooth $k$-varieties we will use $f^*: A^*(Y)\ra A^*(X)$ to indicate both flat and lci pullbacks. Lastly, for a quotient $X/G$ of a smooth variety by the action of a finite group, we write $A^*(X/G)$ for the equivariant Chow ring (which is well defined due to \cite[Theorems 3 and 4]{edidin-graham}).
\subsubsection*{Acknowledgments}
I wish to thank B. Fantechi for suggesting me this problem and for her guidance and support during the preparation of this work. Special thanks to A. Ricolfi for his comments on a previous draft of the paper and for many insightful suggestions. The author is part of INDAM-GNSAGA.
\section{Background}\label{sec-background}
In what follows $X$ is assumed to be a $d$-dimensional, smooth projective variety over a ground field $k=\bar{k}$. For reasons which will be clear later, we assume that 2 and 3 do not divide the characteristic of $k$. Let $\mathcal {Z}_m\subseteq X\times\hilb^m(X)$ be the universal family. As a closed subset it is identified with $$\mathcal{Z}_m:=\{(p,Z)\in X\times\hilb^m(X)\,|\,p\in Z\}.$$
In the case $m=2$ there is an isomorphism 
\begin{equation}\label{equnivf}
    \cal Z_2\cong\bl_\Delta(X\times X),
\end{equation}
which carries the projection onto the first factor to the restriction of the projection $p_X|_{\cal Z_2}:\cal Z_2\subseteq X\times \hilb^2(X)\rightarrow X$. Through this isomorphism, $\hilb^2(X)$ is identified with the quotient $\bl_\Delta(X\times X)/\mathfrak{S}_2$, where the $\mathfrak{S}_2$-action permutes the factors and has the exceptional divisor as fixed locus.

The key to the computation of the ring structure of $A^*(\hilb^3(X))$ relies in observing that the rational morphism $$X\times\hilb^2(X)\dashrightarrow\hilb^3(X)$$
sending the couple $(p,Z)$ to the subscheme $\{p\cup Z\}$ is not defined along the universal family $\cal Z_2\subseteq X\times\hilb^2(X)$. Thus, by blowing up, we obtain the following morphisms
$$
\begin{tikzcd}
    &\bl_{\cal Z_2}(X\times\hilb^2(X))\arrow[dl, "\pi"']\arrow[dr, "\pi_3"]\\
    X\times\hilb^2(X) & &\hilb^3(X).
\end{tikzcd}
$$
For any integer $m\in\Z_{>0}$ define the \textit{nested Hilbert scheme} $\hilb^{[m,m+1]}(X)$ as the closed subset $$\hilb^{[m,m+1]}(X):=\{(Z_m,Z_{m+1})\in\hilb^m(X)\times\hilb^{m+1}(X)\,|\, Z_m\subseteq Z_{m+1}\}.$$
Denote by 
$$\hilb^{m}(X)_{m'}\subseteq\hilb^m(X)$$ the closed subvariety parametrizing subschemes supported at $m'$ points at most (of course $m'\leq m$). Using the same notation as in \cite{Gottsche1993}, there is a residue morphism 
\begin{equation}\label{eqres}
\res:\hilb^{[m,m+1]}(X)\rightarrow X
\end{equation} sending a pair $(Z_m,Z_{m+1})$ to $\operatorname{Supp}(I_{Z_m}/I_{Z_{m+1}})$, i.e. the point in which $Z_{m+1}$ has higher length than $Z_m$. For $X$ smooth projective, so is $\hilb^{[2,3]}(X)$ (see \cite{cheah1998cellular} and \cite{monavari2023lissite}).

We are interested in the case $m=2$. Consider the morphism 
$$\epsilon:=\left(\res, p_{\hilb^2(X)}\right):\hilb^{[2,3]}(X)\rightarrow X\times\hilb^2(X),$$
which lets us define the following subschemes:
\begin{itemize}
    \item $D:=\epsilon^{-1}(X\times\hilb^2(X)_1)$, and 
    \item $F:=\{(Z_2,Z_3)\in\hilb^{[2,3]}(X)\,|\,\res(Z_2,Z_3)\in Z_2\}$, 
\end{itemize}
both being considered with reduced scheme structure.
\begin{prop}[\cite{gottsche2006hilbert}, Proposition 2.5.8]\label{propGbl}
    The map $\epsilon:\hilb^{[2,3]}(X)\rightarrow X\times\hilb^2(X)$ is the blowup along the closed subvariety $\cal Z_2$, and $F\subseteq\hilb^{[2,3]}(X)$ is the exceptional divisor. Through the isomorphism $\hilb^{[2,3]}(X)\cong\bl_{\cal Z_2}(X\times\hilb^2(X))$, the map $\pi_3$ is identified with the projection onto the second factor. 
\end{prop}
To conclude, define 
\begin{itemize}
    \item $D^2:=D\cap F\cong \pi_3^{-1}(\hilb^3(X)_1)$, and 
    \item $\hilb^3(X)^s\subseteq\hilb^3(X)_1$ the closed subvariety parametrizing non curvilinear subschemes, i.e. those whose embedding dimension is strictly greater than one.
\end{itemize}
As shown in the diagram
$$
\begin{tikzcd}
    \hilb^{[2,3]}(X)\arrow[rr,"(\res\text{,}\pi_3)"]\arrow[d]&&X\times\hilb^3(X)\\
    \cal Z_3\arrow[urr]
\end{tikzcd}
$$
the map $(\res,\pi_3):\hilb^{[2,3]}(X)\rightarrow X\times \hilb^3(X)$ factors through the universal family $\cal Z_3\subseteq X\times\hilb^3(X)$. The vertical arrow then is an isomorphism outside of $\pi_3^{-1}(\hilb^3(X)^s)$. As a result, The restriction $$\pi_3|_{\hilb^{[2,3]}(X)\setminus D^2}:\hilb^{[2,3]}(X)\setminus D^2\rightarrow\hilb^3(X)\setminus\hilb^3(X)_1$$
is a a finite and flat cover of degree 3, branched along $\hilb^3(X)_2\setminus\hilb^3(X)_1$. The ramification divisor is $F\setminus D^2$, and the restriction to $D\setminus D^2$ and $F\setminus D^2$ are isomorphisms (see \cite[Section 1]{Gottsche1993}). It follows that the map $$\pi_{3*}\pi_{3}^*:A^*(\hilb^3(X))\rightarrow A^*(\hilb^3(X))$$
equals multiplication by 3. Then, up to taking the coefficients of the Chow groups in $\Z\left[\frac{1}{6}\right]$, we have that $\pi_3^*:A^*(\hilb^3(X))\rightarrow A^*(\hilb^{[2,3]}(X))$ and $\pi_{3*}:A^*(\hilb^{[2,3]}(X))\rightarrow A^*(\hilb^3(X))$ are respectively an injective and surjective homomorphisms of rings. 

The strategy to compute $A^*(\hilb^3(X))$ is then clear. Due to Proposition \ref{propGbl} we obtain a description of $A^*(\hilb^{[2,3]}(X))$ as an algebra with generators and relations over $A^*(X\times\hilb^2(X))$, which is itself an algebra over $A^*(X\times\sym^2(X))$. Then, one can proceed to compute the image of $$\pi_{3*}\pi_3^*:A^*(\hilb^{[2,3]}(X))\rightarrow A^*(\hilb^{[2,3]}(X))$$ on a set of generators. Since we can identify $A^*(\hilb^3(X))$ with its image via $\pi_3^*$ in $A^*(\hilb^{[2,3]}(X))$, we get an explicit presentation of it.
\section{Intermediate step: computing \texorpdfstring{$A^*(\hilb^{[2,3]}(X))$}{A*(Hilb^[2,3](X))}}
\subsection{The Chow ring of a blowup}
We start by recalling the structure of the Chow ring of a blowup of a smooth projective variety $X$ along a smooth closed subvariety $Y$ of codimension $d$. 

Denote by $\iota:Y\hookrightarrow X$ the closed immersion corresponding to $Y$, by $\pi:\bl_Y X\rightarrow X$ the blowup map, and by $i:E\hookrightarrow \bl_Y X$ the inclusion of the exceptional divisor. Lastly, define $p:E\rightarrow Y$ to be the restriction of $\pi$ to the exceptional locus, which is the projective bundle $\pp(N_{Y/X})\ra Y$.

For what follows we will assume $\iota_*:A^*(Y)\rightarrow A^*(X)$ to be injective and $\iota^*:A^*(X)\rightarrow A^*(Y)$ to be surjective.
By the exact sequence for an open immersion \cite[Proposition 1.8]{fulton2013intersection}, we have the following  commutative diagram with exact rows 
\begin{equation}\label{eqopenclosed}
\begin{tikzcd}
    & A^*(E)\arrow[r,"i_*"]\arrow[d, two heads,"p_*"']& A^*(\bl_YX)\arrow[r]\arrow[d,"\pi_*"]& A^*(X\setminus Y)\arrow[r]\arrow[d, equal] & 0\\
    0\arrow[r]&A^*(Y)\arrow[r,"\iota_*"] & A^*(X)\arrow[r] & A^*(X\setminus Y)\arrow[r] & 0.
\end{tikzcd}
\end{equation}
It is well known that the Chow ring of the projective bundle $E=\pp(N_{Y/X})$ is in a natural way an algebra over $A^*(Y)$ with one generator, namely 
\begin{equation}\label{eqnb}
A^*(E)\cong A^*(Y)[h]/(h^d+c_1(N_{Y/X})h^{d-1}+\dots+c_d(N_{Y/X})),
\end{equation}
with $h=c_1(\oo_{\pp(N_{Y/X})}(1))$ \cite[Theorem 3.3]{fulton2013intersection}. Moreover, if $\alpha=\sum_{j=0}^{d-1} \alpha_j \,h^j\in A^*(E)$ (with $\alpha_j\in A^*(Y)$) by \cite[Example 3.3.3]{fulton2013intersection} we have
\begin{equation}\label{eqbl0}
p_*(a\cdot h^j)= \begin{cases} \,0 & \mbox{if } j<d-1\\ \, a & \mbox{if } j=d-1,\end{cases}
\end{equation}
so that $p_*\alpha=\alpha_{d-1}$.
By assumption, the pushforward map 
\begin{equation}
    \iota_*:A^*(Y)\rightarrow A^*(X)
\end{equation}
is injective.
In particular, applying the five-lemma to \eqref{eqopenclosed}, it follows that $\pi_*$ is surjective as well.

Next, we claim that $\ker(i_*)$ is trivial. Indeed, letting $\alpha:=\sum_{j=0}^{d-1}\alpha_j\, h^j$, if $i_*\alpha=0$
then $$0=\pi_* i_*\alpha=\iota_*( p_*\alpha).$$ 
By injectivity of $\iota_*$ and the identity in \eqref{eqbl0}, we get that $0=p_*\alpha=\alpha_{d-1}$. Since $i_*$ is a homomorphism of $i^*(A^*(\bl_Y X))$-modules, and the class $e\in A^*(\bl_Y X))$ of the exceptional divisor pulls back to $-h$, we have that $$i_*(h\cdot\alpha)=i_*(i^*(-e)\cdot\alpha)=-e\cdot (i_*\alpha)=0,$$
which in turn implies $\alpha_{d-2}=0$. Repeating the argument $(d-2)$-times the claim follows. As a result, we have the following commutative diagram with exact rows 
$$
\begin{tikzcd}
    &0\arrow[d]&0\arrow[d]\\
    &\ker  p_*\arrow[d]&\ker\pi_*\arrow[d]\\
    0\arrow[r]&A^*(E)\arrow[r,"i_*"]\arrow[d, two heads,"p_*"']& A^*(\bl_Y X)\arrow[r]\arrow[d, two heads,"\pi_*"]& A^*(X\setminus Y)\arrow[r]\arrow[d, equal] & 0\\
    0\arrow[r]&A^*(Y)\arrow[r,"\iota_*"]\arrow[d] & A^*(X)\arrow[d]\arrow[r] & A^*(X\setminus Y)\arrow[r] & 0\\
    &0&0
\end{tikzcd}
$$
which by the snake-lemma gives an isomorphism $\ker (p_*)\cong\ker(\pi_*)$. The composition $\pi_*\pi^*:A^*(X)\rightarrow A^*(X)$ equals the identity morphism, as for $\alpha\in A^*(X)$ 
$$\pi_*\pi^*(\alpha)=\alpha\cdot\pi_*([\bl_Y X])=\alpha$$
by projection formula. This gives a splitting of the second column 
\begin{equation}\label{eqbl2}
A^*(\bl_Y X)\cong\pi^*(A^*(X))\oplus\bigoplus_{j=0}^{d-2}i_*(A^*(Y)\cdot h^j).
\end{equation}
To conclude, one needs to give a description of the $A^*(X)$-algebra structure on $A^*(\bl_Y X)$. Now, $i_*(A^*(Y))$ is an algebra over $A^*(X)$ through the map $\iota^*: A^*(X)\rightarrow A^*(Y)$, which is assumed to be surjective. Denote $I_{Y^*}:=\ker \iota^*$, and let $\gamma\in I_{Y^*}$. Then 
$$\pi^*(\gamma)\cdot e=\pi^*(\gamma)\cdot i_*(-h)=-i_*(h\cdot p^*(\iota^*\gamma))=0,$$ 
and we have the relation $$\pi^*(I_{Y^*})\cdot \,e=0\quad\hbox{in }A^*(\bl_Y X).$$
From the identity \eqref{eqbl2}, we get 
\begin{equation}\label{eqbl3}
e^d=\pi^*\pi_*(e^d)+\sum_{j=1}^{d-1}\beta_j\cdot e^j=(-1)^{d+1}\pi^*([Y])+\sum_{j=1}^{d-1}\beta_j\cdot e^j,
\end{equation}
for $\beta_1,\dots ,\beta_{d-1}\in A^*(X)$. Following \eqref{eqbl3} and the fact that $$h^d+\sum_{j=0}^{d-1}c_{d-i}(N_{Y/X})\,h^i=0 $$
in $A^*(E)$, we can proceed to describe the coefficients $\beta_i$ as follows. By diagram chasing, and observing that $\pi_*( e^{d+1})=\beta_{d-1}\cdot[Y]$ we obtain $$(-1)^{d+2}\,\beta_{d-1}\cdot [Y]=\pi_*(e^{d+1})=\iota_*p_*i^*((-1)^{d+1}h^{d+1})=\iota_*((-1)^{d+1}c_1(N_{Y/X})),$$
so that $\beta_{d-1}=-c_1(N_{Y/X})$.
Using the same argument for $\pi_*(e^{d+j}),\,j=2,\dots,\,d-1$, we get that 
\begin{equation}\label{eqblx2}
   e^d+\sum_{j=1}^{d-1}c_{d-j}(N_{Y/X})\,e^j+(-1)^d\pi^*[Y]=0,
\end{equation}
in $A^*(\bl_Y X)$. Putting things together we have showed the following.
\begin{prop}\label{propchowbl}
    Let $X$ be a smooth projective variety, and let $\iota:Y\hookrightarrow X$ be the inclusion of a smooth subvariety such that the pushforward $\iota_*:A^*(Y)\rightarrow A^*(X)$ is injective and the pullback $\iota^*:A^*(X)\rightarrow A^*(Y)$ is surjective. Then the Chow ring of $\bl_Y X$ is an algebra over $A^*(X)$ with one generator, the class $e$ of the exceptional divisor, and the relations
    \begin{align}
        &\pi^*(I_{Y^*})\cdot e  = 0, \\
        e^d+\sum_{j=1}^{d-1}c_{d-j}&(N_{Y/X})\,e^j+(-1)^d\pi^*[Y]=0.
    \end{align}
\end{prop}
\subsection{Applications to nested Hilbert schemes}
As anticipated in Section \ref{sec-background}, by \cite[Proposition 2.5.8]{gottsche2006hilbert}, there is an isomorphism $$\bl_{\cal Z_2}(X\times\hilb^2(X))\cong \hilb^{[2,3]}(X).$$
In particular, this makes it feasible to give a concrete description of the Chow ring of the latter. We need a few preliminary lemmata.

Denote 
$$
\begin{tikzcd}
    & X\times\hilb^2(X)\arrow[dl,"p_X"']\arrow[dr, "p_{\hilb^2(X)}"]&\\
    X&&\hilb^2(X)
\end{tikzcd}
$$
the projections to the two factors in the product $X\times\hilb^2(X)$.
\begin{lm}\label{lmz00}
    In $A^*(X\times\hilb^2(X))$ we have 
    \begin{equation}\label{eqcZ}    c\left(N_{\cal Z_2/X\times\hilb^2(X)}^\vee\right)=p_X^*(\,c(T_{X}^\vee))\cdot\frac{1-2e}{1-e},
    \end{equation}
    where $e$ is the class of the exceptional divisor in $A^*(X\times\hilb^2(X))$. 
\end{lm}
\begin{proof}
    Denote $E$ the exceptional divisor in the blowup $\cal Z_2\cong\bl_\Delta(X\times X)$. As in \cite[Lemma 2.1]{Gottsche1993}, it suffices to observe that we have an exact diagram
    $$
    \begin{tikzcd}[column sep=1.5em, row sep=1.5em]
       & & 0\arrow[d] & 0\arrow[d]&\\
         & & p_{\hilb^2(X)}^*(T_{\hilb^2(X)}^\vee)\arrow[d]\arrow[r, equal] &p_{\hilb^2(X)}^*(T_{\hilb^2(X)}^\vee)\arrow[d] & \\
        0\arrow[r] & N_{\cal Z_2/X\times\hilb^2(X)}^\vee\arrow[d, equal]\arrow[r]&(T_{X\times\hilb^2(X)}^\vee)|_{\cal Z_2}\arrow[d]\arrow[r]& T_{\cal Z_2}^\vee\arrow[d]\arrow[r]&0\\
        0\arrow[r] & N_{\cal Z_2/X\times\hilb^2(X)}^\vee\arrow[r] & p_X^*(T_X^\vee)|_{\cal Z_2}\arrow[r]\arrow[d]&\oo_E(-E)\arrow[r]\arrow[d]&0\\
        & & 0 & 0 &
    \end{tikzcd}
    $$
    where the middle row is the sequence for the cotangent bundle of a subvariety and the last column is that of a branched twofold cover. As a result, one may compute $c\left(N_{\cal Z_2/X\times\hilb^2(X)}^\vee\right)$ from the last row, getting \eqref{eqcZ}.
\end{proof}
Let us fix some notation. Denote the product $X^3$ as $X_0\times X_1\times X_2$, write
$$p_i:X_0\times X_1\times X_2\rightarrow X_i$$ 
for the projection to the $i$-th factor for $i=0,1,2$, and likewise for $p_{ij}$. We denote by $\delta_{ij}:\Delta_{ij}\hookrightarrow X_i\times X_j$ the embedding of the diagonal in the product $X_i\times X_j$ and by $\delta_{ijk}:\Delta_{012}\hookrightarrow X_0\times X_1\times X_2$ the embedding of the small diagonal in $X^3$. Therefore, we will identify $X\times\hilb^2(X)$ as  
$$X_0\times (\bl_{\Delta_{12}}(X_1\times X_2)/\mathfrak{S}_2),$$
with the $\mathfrak{S}_2$-action permuting the factors 1 and 2 and fixing the exceptional locus. Lastly, let $$\iota:\cal Z_2\hookrightarrow X\times\hilb^2(X)$$
be the inclusion of the universal family.
\begin{lm}\label{lmz1}
    There is an isomorphism 
    \begin{equation}\label{equniv0}
        \iota_*(A^*(\cal Z_2))\cong (A^*(\bl_{\Delta_{012}}(X_2\times\Delta_{0,1}))\oplus A^*(\bl_{\Delta_{012}}(X_1\times \Delta_{02})))^\mathfrak{S_2}
    \end{equation}
    in $A^*(X\times\hilb^2(X))\cong A^*(X_0\times (\bl_{\Delta_{12}}(X_1\times X_2)/\mathfrak{S}_2))$ with coefficients in $\Z\left[\frac{1}{2}\right]$.
\end{lm}
\begin{proof}
As in \eqref{equnivf}, we can identify $\cal Z_2\cong\bl_\Delta(X\times X)$ and $\hilb^2(X)$ with the quotient $\bl_\Delta(X\times X)/\mathfrak{S}_2$. Observe also that by definition $\cal Z_2\cong\hilb^{[1,2]}(X)$, so that the residue mapping in \eqref{eqres} is defined. As a result, we have two morphisms $\cal Z_2 \rightarrow X\times\hilb^{[1,2]}(X)\cong X_0\times\bl_{\Delta_{12}}(X_1\times X_2)$
\begin{align*}
    j_1:\,(p,Z)&\mapsto(p,p,Z)\\
    j_2:\,(p,Z)&\mapsto(p,\res(p,Z),Z),
\end{align*}
whose image is identified respectively with the two factors appearing in \eqref{equniv0}. In particular, the inclusion of the universal family in $X\times\hilb^2(X)$ for $i=1,2$ factors as follows 
$$
\begin{tikzcd}
    &X\times\hilb^{[1,2]}(X)\arrow[d,"g"]\\
    \cal Z_2\arrow[ur,"j_i"]\arrow[r, "\iota"', hook] & X\times\hilb^2(X),
\end{tikzcd}
$$
where $g=(\id_X\times\varpi):X_0\times\bl_{\Delta_{12}}(X_1\times X_2)\rightarrow X\times\hilb^{2}(X)$ is the product of the identity and the quotient map $\varpi$. Observe that $g^* g_*=\id_{X\times\hilb^{[1,2]}(X)}+\sigma^*$ as an endomorphism in $A^*(X\times\hilb^{[1,2]}(X))$, where $\sigma\in\mathfrak{S}_2$ is the nontrivial element. This means that $\iota=g\circ j_1 =g\circ j_2$, and in particular $$g^*\iota_*=g*(g\circ j_i)_*=g^*g_*j_{i*}=j_i+\sigma^*j_{i*}=j_{1*}+j_{2*}.$$
Finally, we see that $\iota_*=g_* g^*\iota_*=g_*(j_{1*}+j_{2,*})=(\id_X+\varpi)_*(j_{1*}+j_{2*})$, which gives \eqref{equniv0}.
\end{proof}
\begin{thm}\label{chow23}
The Chow ring $A^*(\hilb^{[2,3]}(X))$ with coefficients in $\Z\left[\frac{1}{2}\right]$ is an algebra over $A^*(X\times\sym^2(X))$ with two generators, $e,\,f$, and the following relations
\begin{align}
        I_{\cal Z_2^*}\cdot f  &= 0,\quad I_{\cal Z_2^*}:=\ker((\cal Z_2\hookrightarrow X\times\hilb^2(X)^*)), \label{rel0} \\
        I_{(X\times\Delta_{12})^*}\cdot e & = 0,\quad I_{(X_0\times\Delta_{12})^*}:=\ker((X_0\times\Delta_{12}\hookrightarrow X^3)^*),\label{rel1}\\
        e^d+\sum_{j=1}&^{d-1}c_{d-j}(T_X)\,e^j+(-1)^d[X_0\times\Delta_{12}]=0,\label{rel2}\\
        f^d+\sum_{j=1}&^{d-1}c_{d-j}\left(N_{\cal Z_2/X\times\hilb^2(X)}\right)\,f^j+(-1)^d\left[\frac{1}{2}(X_2\times\Delta_{01}+X_1\times\Delta_{02})\right]=0.\label{rel3}
    \end{align}
\end{thm}
\begin{proof}
    Due to Lemma \ref{lmz1} we cannot apply directly Proposition \ref{propchowbl}, as the kernel of $\iota_*$ may be non zero. In fact, $\iota_*$ kills those classes whose image is anti-symmetric under the involution permuting the factors $1$ and $2$ in the product $X_0\times X_1\times X_2$. 
    However, a similar description may be achieved noticing that $A^*(\pp(N_{\cal Z_2/X\times\hilb^2(X)}))/p^*(\ker\iota_*)$ is still an algebra over $A^*(\cal Z_2)$ with one generator and the relation \eqref{eqnb}. Indeed, it is enough to observe that by Lemma \ref{lmz00} the Chern classes 
    \begin{equation}\label{eqch}
    c_j(N_{\cal Z_2/X\times\hilb^2(X)})\not\in\ker\iota_*\quad\hbox{for each }j=0,\dots d,
    \end{equation} 
    as they are products of pullbacks of classes coming from the $0$-th factor and a polynomial in the class $e$, which is fixed by the $\mathfrak{S}_2$-action. Thus it may be explicitly described as $$A^*(\pp(N_{\cal Z_2/X\times\hilb^2(X)}))/p^*(\ker\iota_*)\cong (A^*(Z_2)/\ker\iota_*)[t]/(t^d+\bar{c}_1t^{d-1}+\dots+\bar{c}_d),$$
    with $\bar{c_j}:=c_j(N_{\cal Z_2/X\times\hilb^2(X)})$. Then, the exact argument as in proposition \ref{propchowbl}, using the identity \eqref{eqbl0}, shows that the map $$A^*(\pp(N_{\cal Z_2/X\times\hilb^2(X)}))/p^*(\ker\iota_*)\rightarrow A^*(\hilb^{[2,3]}(X))$$
    is injective as well. Thus, we have the following diagram
    \begin{center}
    $$
\begin{tikzcd}[column sep=1.1em, row sep=1.5em] 
&0\arrow[d]&0\arrow[d]\\
&\ker \tilde{p}_*\arrow[d]&\ker\pi_*\arrow[d]\\
0\arrow[r]&
\frac{A^*(\pp(N_{\cal Z_2/X\times\hilb^2(X)}))}{p^*(\ker\iota_*)} 
\arrow[r,"i_*"]\arrow[d, two heads,"\tilde{p}_*"']&
A^*(\hilb^{[2,3]}(X))\arrow[r]\arrow[d, two heads,"\pi_*"]&
A^*(\hilb^{[2,3]}(X)\setminus \cal Z_2)\arrow[r]\arrow[d, equal] & 0\\
0\arrow[r]&
A^*(\cal Z_2)/\ker\iota_* 
\arrow[r,"\iota_*"]\arrow[d] &
A^*(X\times\hilb^2(X))\arrow[d]\arrow[r] &
A^*(X\times\hilb^2(X)\setminus \cal Z_2)\arrow[r] & 0\\
&0&0.
\end{tikzcd}
$$
\end{center}
    So, applying the snake-lemma, we conclude that there is an isomorphism between $$\ker(\tilde{p}_*):=\ker(A^*(\pp(N_{\cal Z_2/X\times\hilb^2(X)}))/p^*(\ker\iota_*)\rightarrow A^*(\cal Z_2)/(\ker\iota_*)$$ and $\ker(\pi_*):=\ker(A^*(\hilb^{[2,3]}(X))\rightarrow A^*(X\times\hilb^2(X)))$.
    This yields the decomposition 
    \begin{equation}\label{eqds}
    A^*(\hilb^{[2,3]}(X))\cong A^*(X\times\hilb^2(X))\oplus\bigoplus_{j=1}^{d-1}\,(A^*(\cal Z_2)/\ker\iota_*)\cdot f^j.      
    \end{equation}
    The $A^*(X\times\hilb^2(X))$-algebra structure on $A^*(\cal Z_2)/\ker(\iota_*)$ is obtained by the pull back morphism $\iota^*:X\times\hilb^2(X)\rightarrow A^*(\cal Z_2)$, which is onto due to lemma \ref{lmz1}. As a result, the only relation in this algebra structure is still \eqref{rel0}. To get \eqref{rel3} it is enough to follow the argument in the last part of the proof of proposition \ref{propchowbl}. Indeed, since the observation in \eqref{eqch} grants us that Chern classes of the normal bundle to the universal family $\cal Z_2$ are not annihilated by $\iota_*$, we get an identity as in \eqref{eqblx2}. Lastly, the fact that $\pi^*[\iota_*\cal Z_2]=\left[1/2(X_2\times\Delta_{01}+X_1\times\Delta_{02})\right]$ follows from lemma \ref{lmz1}.

    Looking at the first factor in \eqref{eqds}, using the identity $\bl_\Delta(X\times X)/\mathfrak{S}_2\cong\hilb^2(X)$ there is an isomorphism $$A^*(X\times\hilb^2(X))\cong A^*(\bl_{X_0\times\Delta_{12}}(X_0\times X_1\times X_2))^{\mathfrak{S}_2},$$
    where the $\mathfrak{S}_2$-action permutes the factors 1 and 2. To compute the latter ring, notice that the maps $$X_0\times X_i\cong X_0\times\Delta_{12}\xhookrightarrow{\id\times\delta_{12}} X_0\times X_1\times X_2\xrightarrow{p_{0i}}X_0\times X_i,\quad i=1,2,$$
    compose to the identity, so that by functoriality of the Chow ring
    $$(\id\times\delta_{12})_*:A^*(X_0\times\Delta_{12})\hookrightarrow A^*(X_0\times X_1\times X_2)$$
    is injective, and 
    $$(\id\times\delta_{12})^*:A^*(X_0\times X_1\times X_2)\twoheadrightarrow A^*(X_0\times\Delta_{1,2})$$
    is surjective. Thus, we may apply Proposition \ref{propchowbl}, so that we can further decompose
    \begin{equation}
        A^*(X\times\hilb^2(X))\cong A^*(X\times\sym^2(X))\oplus\bigoplus_{k=1}^{d-1}(\id\times\delta_{12})*(A^*(X_0\times\Delta_{12}))\cdot e^k.
    \end{equation}
    This also grants us automatically the remaining relations. The first one, namely \eqref{rel0}, follows without further work. For \eqref{rel1}, it suffices to observe that $$N_{X_0\times\Delta_{12}/X_0\times X_1\times X_2}\cong p_{12}^*(N_{\Delta_{12}/X_1\times X_2})\cong p_{12}^*(T_X),$$
    and that Chern classes of the latter belong to $A^*(\Delta_{12})$, so that their pullback is invariant under the $\mathfrak{S}_2$-action.
\end{proof}
\section{$A^*(\hilb^3(X))$ for $X$ projective}
The discussion of the previous section allows us to compute the Chow ring of $\hilb^3(X)$. Theorem \ref{chow23} gives a set of generators for $A^*(\hilb^{[2,3]}(X))$; the idea then is to compute the image of such elements under the endomorphism $\pi_3^*\pi_{3*}:A^*(\hilb^{[2,3]}(X))\rightarrow A^*(\hilb^{[2,3]}(X))$ and present $A^*(\hilb^3(X))$ as a subring of $A^*(\hilb^{[2,3]}(X))$. Throughout this section we will make a consistent use of moving lemmas. 
\begin{defi}
    Let $X$ be a nonsingular variety over an algebraically closed field. Two cycles $\alpha,\,\beta$ on $X$ are said to meet \textit{properly} if for each variety $V$ (resp. $W$) which appears with non-zero coefficient in $\alpha$ (resp. $\beta$), $V$ meets $W$ properly, i.e. $\dim(V\cap W)=\dim V+\dim W-\dim X$.
\end{defi}
\begin{thm}[Moving lemma, \cite{fulton2013intersection}, Section 11.4]\label{moving_lemma}
    If $X$ is non-singular and quasi-projective, and $\alpha$ and $\beta$ are cycles on $X$, then there is a cycle $\alpha'$ rationally equivalent to $\alpha$ such that $\alpha'$ meets $\beta$ properly.
\end{thm}
\begin{lm}[\cite{Gottsche1993}, Lemma 3.1]\label{lmpushpull}
    Let $Z\in Z_k(\hilb^3(X))$ be a $k$-cycle.
    \begin{enumerate}
        \item Let $Z':=\pi_3^{-1}(Z\cap\hilb^3(X)_2)\in Z_*(\hilb^{[2,3]}(X))$, and suppose that $\dim(Z')<k$. Then $$\pi_3^*([Z])=\biggl[\overline{\pi_3^{-1}(Z)\setminus Z'}\biggr]\in A^*(\hilb^{[2,3]}(X)).$$
        \item Suppose $Z\subset \hilb^3(X)_1$, and suppose that $Z':=\pi_3^{-1}(Z\cap\hilb^3(X)^s)\in Z_*(\hilb^{[2,3]}(X))$ has dimension strictly less than $k$. Then $$\pi_3^*([Z])=3\biggl[\overline{\pi_3^{-1}(Z)\setminus Z'}\biggr]\in A^*(\hilb^{[2,3]}(X)).$$
        \item Suppose that $Z\subset\hilb^3(X)_2$, and suppose that $Z':=\pi_3^{-1}(Z\cap\hilb^3(X)_1)\in Z_*(\hilb^{[2,3]}(X))$ has dimension strictly less than $k$. Define the cycles $$Z_D:=\overline{(\pi_3^{-1}(Z)\setminus Z')\cap D}\quad\hbox{and}\quad Z_F:=\overline{(\pi_3^{-1}(Z)\setminus Z')\cap F}.$$
        Then $\pi_3^*([Z])=[Z_D]+2[Z_F]$.
    \end{enumerate}
\end{lm}
\subsubsection*{Proof}
    Since both $\hilb^{[2,3]}(X)$ and $\hilb^3(X)$ are smooth and projective we can make use of Theorem \ref{moving_lemma} to prove the claim. The strategy is to check that the cycles have the same intersection numbers with any cycle of complementary dimension. We prove each statement separately.
    \begin{enumerate}
        \item Denote $\hat{Z}:=\overline{\pi_3^{-1}(Z)\setminus Z'}$, and let $W$ be a cycle of complementary dimension. We may assume $W$ to be transversal both to $\hat{Z}$ and $Z'$. Then $\hat{Z}\cap W=\{p_1,\dots,p_r\}$, the intersection multiplicities of $\hat{Z}$ and $W$ at every $p_i$ is $\pm1$, and none of the $p_i$'s lie in $Z'$. We can also assume that $\pi_{3*}[W]=[\pi_3(W)]$, and that the points $\{\pi_3(p_i)\}_i$ are all distinct. By assumption, locally at each $p_i$ the map $\pi_3$ is an isomorphism, so that the multiplicity of $\pi_3(p_i)$ is that of $p_i$ for every $i$. We conclude that $$\pi_{3*}(\pi_3^*[Z]\cdot[W])=\pi_{3*}\biggl[\sum_i m_i p_i\biggr]=\biggl[\sum_i m_i\pi_3(p_i)\biggr]=\pi_{3*}([\hat{Z}]\cdot [W]),$$
        which in turn implies $\pi_3^*([Z])=\hat{Z}$.\\
        \item As before, denote $\hat{Z}:=\overline{\pi_3^{-1}(Z)\setminus Z'}$. Let $W$ be a cycle of complementary dimension which is assumed to intersect transversally $\hat{Z}$ and such that $[\hat{Z}]\cdot[W]=\sum_i m_i[p_i]$. Then the cycle $Z$ can be put within a family $\{Z_t\}_{t\in\mathbb{A}^1}$ such that, locally at each $q_i:=\pi_3(p_i)$, $Z_t$ does not lie in the ramification locus of $\pi_3$. Near each $q_i$ the inverse image $\pi_3^{-1}(Z_t)$ is a union of three connected components, each of which is a deformation of $\hat{Z}$ so that it meets $W$ transversally at points rationally equivalent to the $p_i$'s with multiplicities $m_i$. All of these points can be assumed to have distinct images, thus $[Z_t]\cdot[W]=3\sum_i m_i[p_i]$.\\
        \item The proof is essentially that of the previous item. For a cycle $W$ of complementary dimension, transversal to $Z_D$ and $Z_F$, such that $([Z_D]+[Z_F])\cdot [W]=\sum_i m_i[p_i]$ one can deform the cycle $Z$ in a family $\{Z_t\}_{t\in \mathbb{A}^1}$ such that, locally at $q_i=\pi_3(p_i)$, the cycle $Z_t$ is moved out of the ramification locus. Again, near each $q_i$, the inverse image of $Z_t$ is the union of three connected components, one being a deformation of $Z_D$ and the remaining two a deformation of $Z_F$, meeting $W$ transversally at points with multiplicities $m_i$. We conclude as in the previous point.\qed
    \end{enumerate}
Recall the structure of $A^*(\hilb^{[2,3]}(X))$ as a module over $A^*(X\times\sym^2(X))$ described in Theorem \ref{chow23}. Namely, we have
$$
A^*(\hilb^{[2,3]}(X))\cong A^*(X\times\sym^2(X))\oplus\bigoplus_{k=1}^{d-1}\left((\id\times\delta_{12})_*(A^*(X\times\Delta_{12}))\cdot e^k\right)\oplus\bigoplus_{j=1}^{d-1}\left(A^*(\cal Z_2)/\ker\iota_*\right)\cdot f^j.
$$
Moreover, consider the following commutative diagram 
\begin{equation}\label{eqpi3}
\begin{tikzcd}
    \hilb^{[2,3]}(X)\arrow[d,"\pi_3"']\arrow[r,"\phi"]&X\times\sym^2(X)\arrow[d,"\theta"]\\
    \hilb^3(X)\arrow[r,"\bar{\phi}"]&\sym^3(X),
\end{tikzcd}
\end{equation}
and observe that $\phi^*(A^*(X\times\sym^2(X)))\subset A^*(\hilb^{[2,3]}(X))$ corresponds to the first factor in the previous decomposition.
\begin{prop}\label{prop331}
    We have the following identities in $A^*(\hilb^{[2,3]}(X))$ with coefficients in $\Z\left[\frac{1}{6}\right]$:
    \begin{enumerate}
        \item $\pi_3^*\pi_{3*}(A^*(X\times\sym^2(X)))\cong\phi^*\theta^*(A^*(\sym^3(X));$

                \item $\pi_3^*\pi_{3*}(ef)=3\,ef$;
        \item denote $R:=A^*(\cal Z_2)/(\ker\iota_*)$, then $\pi_3^*\pi_{3*}(R\cdot f)=2(R\cdot f)+2(R\cdot e)$.
        \end{enumerate}
\end{prop}
\subsubsection*{Proof}
    The proof is an application of Lemma \ref{lmpushpull} to each element in the statement.
    \begin{enumerate}
        \item Let $\alpha=[Z]\in A^*(X\times\sym^2(X))$. Observe that the cycle $\pi_3(Z)$ intersects the ramification locus of $\pi_3$ in codimension strictly higher. Therefore, denoting $Z':=\pi_3^{-1}(\pi_3(Z)\cap\hilb^3(X)_2)$, by Lemma \ref{lmpushpull}(1) $$\pi_3^*(\pi_3)_*([Z])=\bigg[\overline{\pi_3^{-1}(\pi_3(Z))\setminus Z'}\bigg].$$
        Next, from \eqref{eqpi3} we have that $\pi_3^*\pi_{3*}\phi^*\,\alpha=\phi^*(\theta^*\theta_*\alpha)$. As a result, $\pi_3^*(\pi_3)_*(A^*(X\times\sym^2(X)))\cong A^*(X\times\sym^2(X))^{\mathfrak{S}_3}$.
        \item The class $ef\in A^*(\hilb^{[2,3]}(X))$ is equal to a multiple of $[\pi_3^{-1}(\hilb^3(X)_1)]$, therefore by Lemma \ref{lmpushpull}(2) we have $$\pi_3^*\pi_{3*}(ef)=3\,ef.$$
        \item Let $\alpha\in R$. Then for the element $\alpha\cdot f=[Z]$ the cycle $\pi_3(Z)$ intersects $\hilb^3(X)_1$ in strictly higher codimension. By Lemma \ref{lmpushpull}(3) it follows that $$\pi_3^*\pi_{3*}([Z])=[Z_D]+2[Z_F]$$
        To conclude it suffices to observe that, by Lemma \ref{lmz1}, $[D]=2\,e$ in $A^*(\hilb^{[2,3]}(X))$.\qed
    \end{enumerate}
In complete analogy with \cite{Gottsche1993}, we have enough data to compute $\pi_3^*\pi_{3*}$ on the full set of generators for $A^*(\hilb^{[2,3]}(X))$ given by Theorem \ref{chow23}.
\begin{prop}\label{proplist}
   Fix $l,\,k\geq 1$, then:
   \begin{itemize}
   \item[(i)] $\pi_3^*\pi_{3*}(a\, e^k f^k)=3\,a\, e^k f^k$ for $a\in A^*(\Delta_{012})$,
   \item[(ii)] $\pi_3^*\pi_{3*}(r\, f^k)=2 \,r\,e^k+2\,r\, f^k+r_0\sum_{i=1}^{k-1}e^i f^{k-i}$ for $r\in R$ and $r_0=\delta_{012}^*(r)$,
   \item[(iii)] $\pi_3^*\pi_{3*}(d\, e^k)=d\, e^k+ d_1\,f^k-d_0\sum_{i=1}^{k-1}e^i f^{k-i}$ for $d\in A^*(X\times\Delta_{12})$, $d_0=\delta_{012}^*(d)$, and $d_1=(\cal Z_2\hookrightarrow X\times\hilb^2(X))^*(d)$,
   \item[(iv)] $\pi_3^*\pi_{3*}(a \,e^k f^{l+k})=a\,(ef)^k(2e^l+2f^l+\sum_{i=1}^{l-1}e^i f^{l-i})$, for $a\in A^*(\Delta_{012})$,
   \item[(v)] $\pi_3^*\pi_{3*}(a\, e^{l+k}f^k)=a\,(ef)^k(e^l+f^l-\sum_{i=1}^{l-1}e^i f^{l-i})$, $a\in A^*(\Delta_{012})$.
   \end{itemize}
\end{prop}
\begin{proof}
    Start by observing that $\pi_3^*\pi_{3*}:A^*(\hilb^{[2,3]}(X))\rightarrow A^*(\hilb^{[2,3]}(X))$ is a homomorphism of $A^*(\hilb^3(X))$-modules. Also, by Proposition \ref{prop331}, the elements $ef$ and $e+f$ belong to $A^*(\hilb^3(X))$, and $\pi_3^*\pi_{3*}(e)=e+f$. Item (i) follows from this observation. Items (iv) and (v) follow from (ii) and (iii) respectively, as $$\pi_3^*\pi_{3*}(a \,e^k f^{l+k})=(ef)^k\pi_3^*\pi_{3*}(a\,f^l)\quad\hbox{and}\quad\pi_3^*\pi_{3*}(a\, e^{l+k}f^k)=(ef)^k\pi_3^*\pi_{3*}(a\, e^l).$$
    To prove (ii) for simplicity set $r=1$. Then $$f^k=(e+f)f^{k-1}-ef^{k-1},$$
    so that $\pi_3^*\pi_{3*}(f^k)=(e+f)\pi_3^*\pi_{3*}(f^{k-1})-ef\pi_3^*\pi_{3*}(f^{k-2})$. To conclude, by induction
    \begin{align*}
        \pi_3^*\pi_{3*}(f^k)&=(e+f)\left(2e^{k-1}+2f^{k-1}+\sum_{i=1}^{k-2}e^if^{k-1-i}\right)-(ef)\left(2e^{k-2}+2f^{k-2}+\sum_{i=1}^{k-3}e^i f^{k-2-i}\right)\\
        &=2e^k+2f^k+\sum_{i=1}^{k-1}e^if^{k-i}.
    \end{align*}
    In the same way, item (iii) follows from the previous inductive argument (here we set $d=1$ again for simplicity):
    \begin{align*}
        \pi_3^*\pi_{3*}(e^k)&=(e+f)\left(2e^{k-1}+2f^{k-1}-\sum_{i=1}^{k-2}e^if^{k-1-i}\right)-(ef)\left(2e^{k-2}+2f^{k-2}-\sum_{i=1}^{k-3}e^i f^{k-2-i}\right)\\
        &=2e^k+2f^k-\sum_{i=1}^{k-1}e^if^{k-i}.
    \end{align*}
    For (ii) and (iii), in the case $r\neq 1$, the same argument applies.
\end{proof}
 Putting everything together, we establish the following.
\begin{thm}\label{thmchow3}
    The Chow ring of $\hilb^3(X),$ with $\Z\left[\frac{1}{6}\right]$-coefficients, is the subring of $A^*(\hilb^{[2,3]}(X))$ generated by those elements of the form 
    \begin{itemize}
        \item $\alpha\in A^*(X\times\sym^2(X))^{\mathfrak{S}_3}$,
        \item $\beta(e^k+f^k)(ef)^m$ with $\beta\in A^*(\Delta_{012})$ and $k\leq (d-1)-m$,
        \item $\gamma(\eta\,e^k+f^k)$ with $\gamma\in R$ and $\eta\in A^*(X\times\Delta_{12})$ and $k\leq d-1$.
    \end{itemize}
\end{thm}
\begin{proof}
    It suffices to notice that by Propositions \ref{prop331} and \ref{proplist} the image of $\pi_3^*\pi_{3*}$ in $A^*(\hilb^{[2,3]}(X))$ is generated by the elements in the statement.
\end{proof}
\begin{rmk}\label{rmkchowkunneth}
    Suppose $X$ admits a Chow--Künneth decomposition, that is $$A^*(X\times X)\cong A^*(X)\otimes A^*(X).$$
    Then we completely recover a result similar to that of \cite{Gottsche1993}. With such hypotheses, the description of $A^*(\hilb^{[2,3]}(X))$ in Theorem \ref{chow23} takes a simpler form. In fact, the latter will be an algebra over $A^*(X)$ with the same two generators $e,\,f$, and the relations \eqref{rel2}, \eqref{rel3}, and $$(1\otimes\xi\otimes1)\cdot e=(1\otimes1\otimes\xi)\cdot e,\quad(\xi\otimes1\otimes1)\cdot f=(1\otimes f\otimes1)\cdot f\quad\hbox{for all }\xi\in A^*(X).$$
    As a result, $A^*(\hilb^3(X))$ can be described, as in \cite[Proposition 3.5]{Gottsche1993}, to be the ring generated by $$(A^*(X)^{\otimes3})^{\mathfrak{S}_3}[e+f,ef]\quad\hbox{and elements of the form}\quad (1\otimes\xi\otimes1)\cdot(e+f),\hbox{ for }\xi\in A^*(X).$$
\end{rmk}
\section{The quasi-projective case}
Knowing the structure of the Chow ring of the Hilbert scheme of three points for smooth projective varieties allows us to compute it also in the quasi-projective case. We will assume the characteristic of $k$ to be zero throughout the section.

Let $X$ be a smooth quasi-projective variety and let $\widetilde{X}$ be a smooth projective compactification (which exists by \cite[\href{https://stacks.math.columbia.edu/tag/0F41}{Tag 0F41}]{stacks-project} and \cite{hironaka1964resolution} since $\operatorname{char}(k)=0$). In particular the inclusion $X\hookrightarrow\wt{X}$ is an open embedding and the square 
\begin{equation}\label{sqopen}
\begin{tikzcd}
    \hilb^m(X)\arrow[r]\arrow[d]&\hilb^m(\wt{X})\arrow[d]\\
    \sym^m(X)\arrow[r]&\sym^m(\wt{X})
\end{tikzcd}    
\end{equation}
is cartesian. Since the arrow on the bottom is an open embedding so is $\hilb^m(X)\hookrightarrow\hilb^m(\wt{X})$. We observe further that the same applies for nested Hilbert schemes, as they can be described in terms of the cartesian square 
$$
\begin{tikzcd}
    \hilb^{[m,m+1]}(X)\arrow[d]\arrow[r]&\hilb^m(X)\times\hilb^{m+1}(X)\arrow[d]\\
    \hilb^{[m,m+1]}(\wt{X})\arrow[r]&\hilb^m(\wt{X})\times\hilb^{m+1}(\wt{X}),
\end{tikzcd}
$$
where the vertical arrow in the right hand side is an open embedding.
Then, the discussion in section \ref{sec-background} for the behaviour of the rational morphism 
$$X\times\hilb^2(X)\dashrightarrow\hilb^3(X)$$
carries verbatim in this context. As a result, we still end up with a generically finite and flat morphism $$\hilb^{[2,3]}(X)\xrightarrow{\pi_3}\hilb^3(X)$$
of degree 3. 
\begin{prop}
    Let $X$ be smooth and quasi-projective. The induced endomorphism at the level of Chow rings 
    $$\pi_{3*}\pi_3^*:A^*(\hilb^3(X))\rightarrow A^*(\hilb^3(X))$$
    equals multiplication by 3.
\end{prop}
\begin{proof}
    We know the statement to hold for smooth, projective varieties. Let $X\hookrightarrow\wt{X}$ be a smooth projective compactification and consider the cartesian square 
    $$
    \begin{tikzcd}
            \hilb^{[2,3]}(X)\arrow[d,"\pi_3"]\arrow[r,"i"]&\hilb^{[2,3]}(\wt{X})\arrow[d,"\wt{\pi}_3"]\\
            \hilb^3(X)\arrow[r,"j"]&\hilb^3(\wt{X}).
    \end{tikzcd}
    $$
    By \eqref{sqopen} the map $j^*:A^*(\hilb^3(\wt{X}))\ra A^*(\hilb^3(X))$ is onto. Let $\eta\in A^*(\hilb^3(X))$ and let $\wt{\eta}\in A^*(\hilb^3(\wt{X})$ be a lift of $\eta$ (i.e. $j^*(\wt{\eta})=\eta$). Then $$\pi_{3*}\pi_3^*(\eta)=\pi_{3*}\pi_3^*(j^*\wt{\eta})=\pi_{3*}i^*\wt{\pi}_3^*(\wt{\eta})=j^*\wt{\pi}_{3*}\wt{\pi}_3^*(\wt{\eta})=j^*(3\,\wt{\eta})=3\,\eta,$$
    so that the statement follows.
\end{proof}
Aa a result, as in the projective case the module and ring structure on $A^*(\hilb^3(X))$ are completely determined by their image in $A^*(\hilb^{[2,3]}(X))$ through $\pi_3^*$. 

In complete analogy with the projective case, Theorem \ref{chow23} holds also for $X$ quasi-projective, as the same proof applies verbatim. So, we still have a ring structure given by 
$$
A^*(\hilb^{[2,3]}(X))\cong A^*(X\times\sym^2(X))\oplus\bigoplus_{k=1}^{d-1}\left((\id\times\delta_{12})_*(A^*(X\times\Delta_{12}))\cdot e^k\right)\oplus\bigoplus_{j=1}^{d-1}\left(A^*(\cal Z_2)/\ker\iota_*\right)\cdot f^j,
$$
together with the relations \eqref{rel0}, \eqref{rel1},\eqref{rel2}, and \eqref{rel3}. 

Consider again the cartesian square
\begin{equation}\label{sqopen1}
    \begin{tikzcd}
            \hilb^{[2,3]}(X)\arrow[d,"\pi_3"]\arrow[r,"i"]&\hilb^{[2,3]}(\wt{X})\arrow[d,"\wt{\pi}_3"]\\
            \hilb^3(X)\arrow[r,"j"]&\hilb^3(\wt{X}),
    \end{tikzcd}
\end{equation}
where $\wt{X}$ is a smooth projective compactification of $X$. 
\begin{rmk}
As before, since $i:\hilb^{[2,3]}(X)\ra\hilb^{[2,3]}(\widetilde {X})$ is an open immersion we have that $i^*:A^*(\hilb^{[2,3]}(\widetilde {X}))\ra A^*(\hilb^{[2,3]}(X))$ is surjective.
\end{rmk}
\begin{lm}\label{lmlift}
    Let $\eta\in A^*(\hilb^{[2,3]}(X))$ and let $\wt{\eta}\in A^*(\hilb^{[2,3]}(\wt{X}))$ be a lift of $\eta$. Then $\pi_3^*\pi_{3*}(\eta)=i^*\wt{\pi}_3^*\wt{\pi}_{3*}(\wt{\eta})$.
\end{lm}
\begin{proof}
    The proof is immediate. Let $\eta$ and $\wt{\eta}$ be as in the statement, then looking at \eqref{sqopen1} we have that $
    \pi_3^*\pi_{3*}(\eta)=\pi_3^*\pi_{3*}i^*(\wt{\eta})=\wt{\pi}_{3*}j^*\pi_3^*(\wt{\eta})=i^*\wt{\pi}_3^*\wt{\pi}_{3*}(\wt{\eta}).$
\end{proof}    
\begin{prop}
    Once again, the following identities hold in $A^*(\hilb^{[2,3]}(X))$:
    \begin{enumerate}
        \item $\pi_3^*\pi_{3*}(A^*(X\times\sym^2(X)))\cong\phi^*\theta^*(A^*(\sym^3(X));$

                \item $\pi_3^*\pi_{3*}(ef)=3\,ef$;
        \item if $R:=A^*(\cal Z_2)/(\ker\iota_*)$, then $\pi_3^*\pi_{3*}(R\cdot f)=2(R\cdot f)+2(R\cdot e)$.
        \end{enumerate}
\end{prop}
\begin{proof}
    By Proposition \ref{prop331} we know the statement to hold for smooth projective varieties. Let $\wt{X}$ be a smooth projective compactification of $X$. For item (1), consider the commutative cube
        \begin{equation}\label{eqcube}
        \begin{tikzcd}[row sep=scriptsize, column sep=scriptsize]
& \hilb^{[2,3]}(X) \arrow[dl, "i"'] \arrow[rr,"\phi"] \arrow[dd] & & X\times\sym^2(X) \arrow[dl, "i'"] \arrow[dd, "\theta"] \\
\hilb^{[2,3]}(\wt{X}) \arrow[rr, "\wt{\phi}"{near end, yshift=1pt}, crossing over] \arrow[dd] & & \wt{X}\times\sym^2(\wt{X}) \\
& \hilb^3(X) \arrow[dl, "j"'] \arrow[rr] & & \sym^3(X) \arrow[dl, "j'"] \\
\hilb^3(\wt{X}) \arrow[rr] & & \sym^3(\wt{X}). \arrow[from=uu, "\wt{\theta}"{near start, xshift=1pt}, crossing over]\\
\end{tikzcd}
        \end{equation}
        By Lemma \ref{lmlift} we know that $(\pi_3^*\pi_{3*})(\phi^*A^*(X\times\sym^2(X)))=i^*(\wt{\pi}_3^*\wt{\pi}_{3*})(\wt{\phi}^*A^*(\wt{X}\times\sym^2(\wt{X})))$. Next, by Proposition \ref{prop331}, $(\wt{\pi}_3^*\wt{\pi}_{3*})(\wt{\phi}^*A^*(\wt{X}\times\sym^2(\wt{X})))=\wt{\theta}^*\wt{\phi}^*(A^*(\sym^3(\wt{X})))$. We conclude by commutativity of the diagram \eqref{eqcube} and surjectivity of $j'^*:A^*(\sym^3(\wt{X})\twoheadrightarrow A^*(\sym^3(X))$ that $$i^*\wt{\theta}^*\wt{\phi}^*(A^*(\sym^3(\wt{X})))=\phi^*\theta^*j'^{*}A^*(\sym^3(\wt{X}))=\phi^*\theta^*A^*(\sym^3(X)).$$
        Items (2) and (3) follow similarly by observing that the classes $e,\,f,\,ef\in A^*(\hilb^{[2,3]})(X)$ admit lifts by the same classes $\wt{e},\,\wt{f},\,\wt{e}\wt{f}\in A^*(\hilb^{[2,3]}(\wt{X}))$. We conclude by Lemma \ref{lmlift}.
\end{proof}
\begin{prop}\label{proplist1}
   Fix $l,\,k\geq 1$, then:
   \begin{itemize}
   \item[(i)] $\pi_3^*\pi_{3*}(a\, e^k f^k)=3\,a\, e^k f^k$ for $a\in A^*(\Delta_{012})$,
   \item[(ii)] $\pi_3^*\pi_{3*}(r\, f^k)=2 \,r\,e^k+2\,r\, f^k+r_0\sum_{i=1}^{k-1}e^i f^{k-i}$ for $r\in R$ and $r_0=\delta_{012}^*(r)$,
   \item[(iii)] $\pi_3^*\pi_{3*}(d\, e^k)=d\, e^k+ d_1\,f^k-d_0\sum_{i=1}^{k-1}e^i f^{k-i}$ for $d\in A^*(X\times\Delta_{12})$, $d_0=\delta_{012}^*(d)$, and $d_1=(\cal Z_2\hookrightarrow X\times\hilb^2(X))^*(d)$,
   \item[(iv)] $\pi_3^*\pi_{3*}(a \,e^k f^{l+k})=a\,(ef)^k(2e^l+2f^l+\sum_{i=1}^{l-1}e^i f^{l-i})$, for $a\in A^*(\Delta_{012})$,
   \item[(v)] $\pi_3^*\pi_{3*}(a\, e^{l+k}f^k)=a\,(ef)^k(e^l+f^l-\sum_{i=1}^{l-1}e^i f^{l-i})$, $a\in A^*(\Delta_{012})$.
   \end{itemize}
\end{prop}
\begin{proof}
    The proof is identical to that of Proposition \ref{proplist}.
\end{proof}
As in the projective case, the Chow ring $A^*(\hilb^3(X))$ admits a description as follows.
\begin{cor}
    Let $X$ be a smooth quasi-projective variety. The Chow ring of $\hilb^3(X),$ with $\Z\left[\frac{1}{6}\right]$-coefficients, is the subring of $A^*(\hilb^{[2,3]}(X))$ generated by those elements of the form 
    \begin{itemize}
        \item $\alpha\in A^*(X\times\sym^2(X))^{\mathfrak{S}_3}$,
        \item $\beta(e^k+f^k)(ef)^m$ with $\beta\in A^*(\Delta_{0,1,2})$ and $k\leq (d-1)-m$,
        \item $\gamma(\eta\,e^k+f^k)$ with $\gamma\in R$ and $\eta\in A^*(X\times\Delta_{1,2})$ and $k\leq d-1$.
    \end{itemize}
\end{cor}
\begin{rmk}
    Once again, if a smooth quasi-projective variety $X$ admits a Chow--Künneth decomposition, we get the same description of $A^*(\hilb^3X))$ as in Remark \ref{rmkchowkunneth}. Observe that this is slightly stronger than the result in \cite{Gottsche1993}. A key feature in the work by Fantechi--Göttsche is that the vector spaces $H^i(X,\mathbb{Q})$ are assumed to be finite dimensional (clearly implied by the projectivity of $X$). If $X$ is only quasi-projective, this may easily fail for $A_i(X)\otimes\mathbb{Q}$. The resulting description of $A^*(\hilb^3(X))$ is less explicit but more intrinsic.
\end{rmk}
\bibliographystyle{alpha}
\bibliography{bib}

@book{fulton2013intersection,
  title={Intersection theory},
  author={Fulton, William},
  volume={2},
  year={2013},
  publisher={Springer Science \& Business Media}
}

@article{edidin-graham,
    AUTHOR = "Edidin, Dan and Graham, William",
    TITLE = "Equivariant intersection theory",
    JOURNAL = "Invent. Math.",
    FJOURNAL = "Inventiones Mathematicae",
    VOLUME = "131",
    YEAR = "1998",
    NUMBER = "3",
    PAGES = "595--634"
}

@article{Gottsche1993,
author = {Fantechi, Barbara and Göttsche, Lothar},
journal = {Journal für die reine und angewandte Mathematik},
pages = {147-158},
title = {The cohomology ring of the {H}ilbert scheme of 3 points on a smooth projective variety.},
url = {http://eudml.org/doc/153527},
volume = {439},
year = {1993},
}

@book{gottsche2006hilbert,
  title={Hilbert schemes of zero-dimensional subschemes of smooth varieties},
  author={G{\"o}ttsche, Lothar},
  year={2006},
  publisher={Springer}
}

@article{grothendieck1957fondements,
  title={Fondements de la g{\'e}om{\'e}trie alg{\'e}brique. Commentaires},
  author={Grothendieck, Alexander},
  journal={S{\'e}minaire Bourbaki},
  volume={7},
  number={297-298},
  pages={13},
  year={1957},
  publisher={Soci{\'e}t{\'e} Math{\'e}matique de France}
}

@misc{stacks-project,
    shorthand    = {Stacks},
    author       = {The {Stacks Project Authors}},
    title        = {\textit{Stacks Project}},
    howpublished = {\url{https://stacks.math.columbia.edu}},
    year         = {2025},
  }

@article{hironaka1964resolution,
  title={Resolution of singularities of an algebraic variety over a field of characteristic zero: II},
  author={Hironaka, Heisuke},
  journal={Annals of Mathematics},
  volume={79},
  number={2},
  pages={205--326},
  year={1964},
  publisher={JSTOR}
}

@article{fogarty1968algebraic,
  title={Algebraic families on an algebraic surface},
  author={Fogarty, John},
  journal={American Journal of Mathematics},
  volume={90},
  number={2},
  pages={511--521},
  year={1968},
  publisher={JSTOR}
}

@article{jelisiejew2020pathologies,
  title={Pathologies on the {H}ilbert scheme of points},
  author={Jelisiejew, Joachim},
  journal={Inventiones mathematicae},
  volume={220},
  number={2},
  pages={581--610},
  year={2020},
  publisher={Springer}
}

@article{ellingsrud1987homology,
  title={On the homology of the {H}ilbert scheme of points in the plane},
  author={Ellingsrud, Geir and Str{\o}mme, Stein Arild},
  journal={Inventiones mathematicae},
  volume={87},
  number={2},
  pages={343--352},
  year={1987},
  publisher={Springer-Verlag Berlin/Heidelberg}
}

@article{gottsche1990betti,
  title={The {B}etti numbers of the {H}ilbert scheme of points on a smooth projective surface},
  author={G{\"o}ttsche, Lothar},
  journal={Mathematische Annalen},
  volume={286},
  number={1},
  pages={193--207},
  year={1990},
  publisher={Springer}
}

@article{DECATALDO2002824,
title = {The {C}how {G}roups and the {M}otive of the {H}ilbert Scheme of {P}oints on a {S}urface},
journal = {Journal of Algebra},
volume = {251},
number = {2},
pages = {824-848},
year = {2002},
issn = {0021-8693},
doi = {https://doi.org/10.1006/jabr.2001.9105},
url = {https://www.sciencedirect.com/science/article/pii/S0021869301991057},
author = {Mark Andrea A {de Cataldo} and Luca Migliorini}
}

@article{nakajima1997heisenberg,
  title={Heisenberg algebra and {H}ilbert schemes of points on projective surfaces},
  author={Nakajima, Hiraku},
  journal={Annals of mathematics},
  volume={145},
  number={2},
  pages={379--388},
  year={1997},
  publisher={JSTOR}
}

@article{graffeo2024motive,
  title={The motive of the {H}ilbert scheme of points in all dimensions},
  author={Graffeo, Michele and Monavari, Sergej and Moschetti, Riccardo and Ricolfi, Andrea T},
  journal={arXiv preprint arXiv:2406.14321},
  year={2024}
}

@article{zhan2022punctual,
  title={Punctual {H}ilbert schemes of points of $\mathbb{A}^3
  $ in the Grothendieck group of varieties},
  author={Zhan, Sailun},
  journal={arXiv preprint arXiv:2208.02419},
  year={2022}
}

@incollection{katz1994desingularization,
  title={The desingularization of $\hilb^4(\mathbb{P}^3)$ and its {B}etti numbers},
  author={Katz, Sheldon},
  booktitle={Zero-dimensional schemes (Ravello, 1992)},
  pages={231--242},
  year={1994},
  publisher={de Gruyter, Berlin}
}

@article{shen2016motive,
  title={The motive of the {H}ilbert cube},
  author={Shen, Mingmin and Vial, Charles},
  journal={Forum of Mathematics, Sigma},
  volume={4},
  pages={e30},
  year={2016},
  organization={Cambridge University Press}
}

@article{cheah1998cellular,
  title={Cellular decompositions for nested {H}ilbert schemes of points},
  author={Cheah, Jan},
  journal={Pacific Journal of Mathematics},
  volume={183},
  number={1},
  pages={39--90},
  year={1998},
  publisher={Mathematical Sciences Publishers}
}

@article{monavari2023lissite,
  title={Sur la lissit{\'e} du sch{\'e}ma {Q}uot ponctuel embo{\^\i}t{\'e}},
  author={Monavari, Sergej and Ricolfi, Andrea T},
  journal={Canadian Mathematical Bulletin},
  volume={66},
  number={1},
  pages={178--184},
  year={2023},
  publisher={Canadian Mathematical Society}
}
\end{document}